\documentclass[centertags]{article}

\usepackage{amsmath}
\usepackage{amssymb}
\usepackage{latexsym}
\usepackage{amsxtra}
\usepackage{amscd}
\usepackage{theorem}
\usepackage{xypic}
\usepackage{epic}
\usepackage{eepic}

\theoremstyle{plain}

{\theorembodyfont{\slshape}

        \newtheorem{thm}{Theorem}[section]
        \newtheorem{cor}[thm]{Corollary}
        \newtheorem{lem}[thm]{Lemma}
        \newtheorem{prop}[thm]{Proposition}

}

{\theorembodyfont{\rmfamily}

        \newtheorem{defn}[thm]{Definition}
        \newtheorem{prob}[thm]{Problem}

        \newtheorem{exa}[thm]{Example}

}

\renewcommand{\em}{\sl}

\newcommand{\proof}{{\bf Proof:\ }}

\newcommand{\Endproof}{\hspace*{\fill} $\Box$ \vspace{1ex} \noindent }

\makeatletter
\renewcommand{\subsection}{\@startsection{subsection}{2}%
        {\z@}{-3.25ex plus -1ex minus-.2ex}{-1em}{\bf}}
\makeatother

\newcommand{\PP}{\mathbb{P}}
\newcommand{\ZZ}{\mathbb{Z}}
\newcommand{\CC}{\mathbb{C}}

\newcommand{\QQ}{\mathbb{Q}}

\newcommand{\FF}{\mathbb{F}}
\renewcommand{\AA}{\mathbb{A}}

\newcommand{\diff}{{\rm d}}

\renewcommand{\a}{{\mathbf a}}
\newcommand{\A}{{\mathbf A}}
\newcommand{\QQb}{{\bar{\QQ}}}

\newcommand{\Aut}{\mathop{\rm Aut}\nolimits}

\newcommand{\ord}{{\rm ord}}

\newcommand{\Res}{\mathop{\rm Res}\nolimits}
\newcommand{\inj}{\hookrightarrow}

\newcommand{\iso}{\stackrel{\sim}{\to}}

\newcommand{\gen}[1]{\mathopen\langle#1\mathclose\rangle}

\newcommand{\ssigma}{{\mathbf \sigma}}

\title{Deformation data, Belyi maps, and the local lifting problem}
\author{Irene I.\ Bouw, Stefan Wewers and Leonardo Zapponi}
\date{}
\begin{document}
\maketitle

\begin{abstract} We prove existence and nonexistence results for certain
  differential forms in positive characteristic, called {\em good deformation
    data}. Some of these results are obtained by reduction modulo $p$ of Belyi
  maps.  As an application, we solve the local lifting problem for groups with
  Sylow $p$-subgroup of order $p$.

\end{abstract}

\section*{Introduction}

\subsection*{The local lifting problem}

Let $k$ be an algebraically closed field of characteristic $p>0$, and
let $G$ be a finite group. A {\em local $G$-action} is a continuous,
faithful and $k$-linear action $\phi:G\inj\Aut_k(k[[x]])$ of $G$ on
the ring of formal power series $k[[x]]$. If such an action exists,
then $G=P\rtimes C$, where $P$ is a $p$-group and $C\cong\ZZ/m\ZZ$ is
a cyclic group of order $m$, with $(m,p)=1$. 

The {\em local lifting problem} concerns the following question. Given
a local $G$-action $\phi$, does there exist a lift of $\phi$ to an
$R$-linear action $\phi_R:G\inj R[[x]]$, where $R$ is a complete
discrete valuation ring of characteristic zero with residue field $k$?
If such a lift exists, we will say $\phi$ {\em lifts to characteristic
  zero}. See \cite{GM2} or \cite{CGH} for a general discussion
of the local lifting problem.

In the previous paper \cite{Dp} we have have shown that, for $G$ the dihedral
group of order $2p$, every local $G$-action can be lifted to characteristic
zero. The methods used in \cite{Dp} apply more generally to the case where the
Sylow $p$-subgroup $P$ of $G$ has order $p$. In this more general situation not
every local $G$-action lifts.  It is easy to write down certain necessary
conditions (\cite{Dp}, Proposition 1.3). We already suspected that these
conditions are also sufficient (\cite{Dp}, Question 1.4) and reduced the proof
of sufficiency to the existence of a certain differential form in
characteristic $p$. This differential form is the deformation datum referred to
in the title of the present paper.

Our first new result says that the differential form needed to prove
sufficiency always exists. As a consequence, the local lifting problem has now
been solved for all groups $G$ with a Sylow $p$-subgroup $P$ of order $p$.

\subsection*{Good deformation data}

Let $G=P\rtimes C$, where $P$ is cyclic of order $p$ and $C$ is cyclic of
order $m$ and where $m$ is prime to $p$. Choose a generator $\tau$ of $P$ and
denote by $\chi:C\to\FF_p^\times$ the character such that
$\sigma\tau\sigma^{-1}=\tau^{\chi(\sigma)}$, for $\sigma\in C$.

Let $\phi:G\inj\Aut_k(k[[x]])$ be a local $G$-action. The {\em
  conductor} of $\phi$ is defined as the (positive) integer
\begin{equation} \label{conductoreq}
  h := \ord_x\left(\frac{\phi(\tau)(x)}{x}-1\right),
\end{equation}
where $\ord_x(f)$ denotes the order of zero of a power series $f\in
k[[x]]$.  It is easy to see that $h$ is prime to $p$. We say that the
local $G$-action $\phi$ has {\em large conductor} if $h>p$ and that it
has {\em small conductor} if $h<p$.

The results of the present paper concern mainly the case of small
conductor. Hence we suppose, for the rest of this introduction, that
$h<p$. 

Let $\phi_R:G\inj\Aut_R(R[[x]])$ be a lift of $\phi$ to characteristic
zero. By a construction going back to \cite{GM} (see also \cite{Henrio} and
\cite{Dp}), we can associate to $\phi_R$ a differential form
$\omega=f(z)\,\diff z$ over a rational function field $k(z)$ with the following
properties:
\begin{itemize}
\item 
  $\omega$ is logarithmic, i.e.\ of the form $\omega=\diff g/g$,
  with $g\in k(z)$,
\item
  $\omega$ has a single zero of order $h-1$.
\item
  There is a faithful action $C\inj\Aut_k(k(z))$ such that 
  \[
         \sigma^*\omega =\chi(\sigma)\cdot\omega,
  \]
  for $\sigma\in C$. 
\end{itemize}
A differential form $\omega$ with these properties is called a {\em good
  deformation datum} with {\em conductor $h$}. (Note that this definition
depends on the character $\chi$, which we regard as fixed.)

Conversely, a good deformation datum $\omega$ as above essentially
determines a lift of $\phi$ to characteristic zero. In particular, a
local $G$-action $\phi$ lifts to characteristic zero if and only if
there exists a good deformation datum with the same conductor $h$.

In Section \ref{sec1} we give necessary conditions for the existence
of a good deformation datum with given conductor (Proposition
\ref{prop0}). Then we show that these conditions are also sufficient,
at least for small conductor (Proposition \ref{prop1}). In Section
\ref{sec2} we show that these two results give a complete solution to
the lifting problem for groups of the form $G=P\rtimes C$, with
$|P|=p$.

\subsection*{Belyi morphisms}

Let us assume, for simplicity, that $m=1$ in the situation considered
above.  Let $\omega=f(z)\,\diff z$ be a good deformation datum with
conductor $h$.  After a change of coordinate, we may assume that the
unique zero of $\omega$ is located at $z=\infty$. Let
$z=z_1,\ldots,z=z_r$ be the positions of the poles of $\omega$. Since
$\omega$ is logarithmic, it has at most simple poles.  Therefore,
$r=h+1$. Moreover, the residue
\[
       a_i := \Res_{z=z_i}(\omega)
\]
is an element of $\FF_p^\times$. Note that $\sum_ia_i=0$ by the residue
theorem. We call the tuple $\a=(a_1,\ldots,a_r)$ the {\em type} of $\omega$.

In Sections 3 to 6 we consider the following problem: given a tuple
$\a=(a_1,\ldots,a_r)\in(\FF_p^\times)^r$ with $\sum_ia_i=0$, does there exists
a good deformation datum $\omega$ of type $\a$? This is a refined version of
the problem discussed above, where we only fix the conductor $h$.

Given $\a=(a_1,\ldots,a_r)$ as above, it is easy to write down a system of
algebraic equations whose solutions correspond to good deformation data of
type $\a$. So far, most attempts to construct good deformation data consisted
in trying to solve these equations. For instance, \cite{Dp}, Theorem 5.6,
relies on the mysterious fact that the system of equations corresponding to a
particular type $\a$ has a unique solution. The proof of this fact given in
\cite{Dp} is quite complicated and offers no explanation why this type $\a$
works and others don't.

In this paper we follow a different path. Consider a good deformation datum
$\omega$ of type $\a=(a_1,\ldots,a_r)$. Then to every choice of lifts
$A_i\in\ZZ$ of $a_i$ such that $\sum_i A_i=0$ one can associate an essentially
unique rational function $g\in k(z)$ with the following properties:
\begin{itemize}
\item
  $\omega=\diff g/g$,
\item
  $g$ has exactly $r$ zeroes/poles, of order $A_1,\ldots,A_r$, and
\item
  $g:\PP^1_k\to\PP^1_k$ is branched exactly over $0$, $1$ and $\infty$. The
  fiber $g^{-1}(1)$ contains a unique ramification point, with ramification
  index $h=r-1$. 
\end{itemize}
If, moreover, $\omega$ has small conductor $h<p$ then $g$ is at most tamely
ramified. It follows from a theorem of Grothendieck that $g$ is the reduction
modulo $p$ of a {\em Belyi map}, i.e.\ a rational function
$f:\PP^1_\CC\to\PP^1_\CC$ with exactly three critical values $0$, $1$,
$\infty$. The ramification structure of $f$ is the same as that of $g$; in
particular, it is completely determined by the tuple $\A=(A_1,\ldots,A_r)$. We
say that $f$ is a Belyi map {\em of type $\A$}.

Our new method consists in choosing a suitable lift $\A$ of $\a$ and to
construct a Belyi map $f$ of type $\A$. One then obtains a good deformation
datum $\omega=\diff g/g$ from the reduction $g$ of $f$ modulo $p$. Whether
this is possible or not for a particular type $\a$ can be easily decided,
thanks to the following lemma.

\vspace{2ex}\noindent
{\bf Lemma}\hspace{0.5em} {\em
  Given $\A=(A_1,\ldots,A_r)\in\ZZ^r$ with $\sum_i A_i=0$; set
  \[
       n_\A:=\sum_i \max(A_i,0), \qquad k_\A={\rm gcd}(A_1,\ldots,A_r).
  \]
  Then there exists a Belyi map of type $\A$ if and only if 
  \[
       k_\A\cdot(r-1) \leq n_\A.
  \]
\vspace{1ex}}

It seems that this lemma has been discovered independently, in various forms,
by several authors. For instance, it is a special case (the {\em genus zero
  case}) of a result of Boccara \cite{Boccara}. It also follows from Theorem 1
of \cite{Zannier}. In Section \ref{belyi} we shall give a simple proof of the
above lemma which uses {\em dessins d'enfants}. This proof is very similar to
the proof given in \cite{Zannier}. We have decided to keep it in the paper,
because our use of dessins may be considered as more geometric and easier to
understand than explicit manipulations with permutations, as in
\cite{Zannier}.

The lemma, combined with Grothendieck's theory of lifting and reduction of
tamely ramified coverings, immediately gives necessary and sufficient
conditions for the existence of good deformation data of a given type $\a$, in
the case where $m=1$. See Corollary \ref{cor2} and Theorem \ref{thm1}.

One can adapt the method described above to the case $m=2$ and hence obtain
almost complete results as well. This is done in Section \ref{m=2}. 

It is not clear how to adapt the method to the case $m>2$. However, one can
still use the case $m=1$ to obtain nonexistence results for good deformation
data with arbitrary $m$. It seems that all nonexistence results that have been
documented in the literature so far can be explained in this way.

\section{Good deformation data} \label{sec1}

Let $k$ be an algebraically closed field of
characteristic $p>0$. We let $K=k(z)$ denote the rational function field over
$k$, which we regard as the function field of the projective line
$X=\PP^1_k$. 

Let $m\geq 1$ be a positive integer which is prime to $p$. We choose, once and
for all, a primitive $m$th root of unity $\zeta_m\in k^\times$. We denote by
$\sigma:X\iso X$ the automorphism of order $m$ such that
$\sigma^*z=\zeta_mz$. 

\begin{defn} \label{def1}
  A {\em deformation datum} is a differential form $\omega=f(z)\,\diff
  z\in\Omega^1_{K/k}$ with the following properties:
  \begin{itemize}
  \item[(a)]
    $\omega$ is an eigenvector for $\sigma$, i.e.\
    \[
        \sigma^*\omega =\zeta_m^c\omega,
    \]
    for some $c\in\ZZ$, and
  \item[(b)]
    $\omega$ is logarithmic, i.e.\ $\omega=\diff g/g$, for some $g\in
    K^\times$.
  \end{itemize}
  A deformation datum $\omega$ is called {\em primitive} if the
  integer $c$ in (a) is prime to $m$. It is called {\em good} if it
  has a unique zero, i.e.\ there exists a unique point $x\in X$ with
  $\ord_x\omega>0$. If this is the case then the integer
  \[
         h:=\ord_x\omega +1
  \]
  is called the {\em conductor} of $\omega$. 
\end{defn}

We remark that in our previous papers the definition of the term
`deformation datum' is more general. Deformation data occur naturally
in a variety of contexts, see e.g.\ \cite{indi}, \cite{habil}. Here we
focus on applications to the local lifting problem, as explained in
the introduction. To save paper we have chosen a more restricted setup
from the start, suitable for the applications we have in mind.

The present paper is concerned with the following problem.

\begin{prob} \label{prob1}
  For which $m$ and $h$ does there exist a good deformation datum $\omega$ of
  conductor $h$?
\end{prob}

The following proposition 
gives necessary conditions on $h$ and $m$. 

\begin{prop} \label{prop0} Fix positive integers $m$ and $h$ and
  assume that there exists a good deformation datum $\omega$ of
  conductor $h$. Then the following holds.
  \begin{enumerate}
  \item
    The conductor $h$ is prime to $p$.
  \item 
    $\omega$ is primitive if and only if $h$ is prime to $m$. If
    this is the case, then
    \[
       m|(p-1)
    \]
    (and hence $\zeta_m\in\FF_p^\times$) and 
    \[
       h \equiv -1 \pmod{m}.
    \]
  \item If $\omega$ is not primitive, then $m|h$ and
    $\sigma^*\omega=\omega$ (i.e.\ $c\equiv 0\pmod{m}$ in Definition
    \ref{def1} (a)).
  \end{enumerate}
\end{prop}

This is a special case of \cite{Dp}, Lemma 3.3.(v). For convenience,
we recall the proof.

\proof Let $x\in X$ be the unique zero of $\omega$. Choose a local
coordinate $w\in k(z)$ at $x$ and a function $g\in k(z)$ such that
$\omega=\diff g/g$. By multiplying $g$ with a $p$th power, if
necessary, we can achieve that $g(x)=1$. Writing $g$ as a power series
in $w$ and computing $\omega=\diff g/g$, one sees that
$h=\ord_x\omega+1\not\equiv 0 \pmod{p}$. This proves (i).

The statements (ii) and (iii) are trivial for $m=1$. We may therefore
assume $m>1$. Since the automorphism $\sigma:X\iso X$ has order $m$,
it has exactly two fixed points, namely $z=0$ and $z=\infty$. The unique
zero $x$ of $\omega$ is clearly fixed by $\sigma$. Replacing the
coordinate $z$ by $z^{-1}$, if necessary, we may assume that $x$ is
the point $z=\infty$ and that $w=z^{-1}$. Now
$\sigma^*w=\zeta_m^{-1}w$, and Condition (a) of Definition \ref{def1}
implies
\begin{equation} \label{prop0eq1}
    h=\ord_x(\omega)+1\equiv -c \pmod{m}.
\end{equation}
In particular, $\omega$ is primitive if and only if $h$ is prime to $m$. 

The same argument used to prove \eqref{prop0eq1} shows that
\begin{equation} \label{prop0eq2}
  \ord_{z=0}(\omega)+1 \equiv c \pmod{m}.
\end{equation}
But $\ord_{z=0}(\omega)$ is either equal to $-1$ or to $0$. In the
first case, $c$ and $h$ are divisible by $m$; this corresponds to Part
(iii) of the proposition. In the second case, $c\equiv 1 \pmod{m}$,
$h\equiv -1\pmod{m}$ and $\omega$ is primitive. This corresponds to
Part (ii) of the proposition.

It remains to prove that $m|(p-1)$ in the second case. Let $z=z_1$ be
a pole of $\omega$ and set $a_1:=\Res_{z=z_1}(\omega)$.  Set
$z_2:=\sigma(z_1)=\zeta_m z_1$.  Condition (a) of Definition
\ref{def1}, together with the congruence $c\equiv 1 \pmod{m}$, implies
that
\[
   a_2:=\Res_{z=z_2}(\omega)=\zeta_m^{-1}a_1.
\]
Since $\omega$ is logarithmic, the residues $a_1,a_2$ actually lie in
$\FF_p^\times\subset k^\times$ (see Lemma \ref{lem1} below).
Therefore, $\zeta_m\in\FF_p^\times$, which is equivalent to
$m|(p-1)$. This finishes the proof of the proposition.
\Endproof

The next result says that the necessary conditions given by
Proposition \ref{prop0} are also sufficient, at least if $h<p$. To
keep the statement simple, we first deal with the case $m|h$ (the
non-primitive case). Here one can immediately write down a good
deformation datum of conductor $h$ (see also \cite{Henrio}, \S 3.5):
\begin{equation} \label{npomega}
    \omega:= \frac{h\,\diff z}{z^{h+1}-z}=\frac{\diff g}{g},
     \qquad \text{with $g:=(z^h-1)/z^h$.}
\end{equation}
It therefore suffices to consider the primitive case (Part (ii) of
Proposition \ref{prop0}). 

\begin{prop} \label{prop1}
  Assume $m|(p-1)$ and let $h$ be a positive integer with 
  \[
      h<p \qquad\text{and}\qquad h \equiv -1 \pmod{m}.
  \]
  Then there exists a good deformation datum with conductor $h$.
\end{prop}

The proof of the proposition is based on the following well-known lemma.

\begin{lem} \label{lem1}
  Let $\omega=f(z)\,\diff z\in\Omega^1_{K/k}$ be a meromorphic differential form
  on $X=\PP^1_k$. Then $\omega$ is logarithmic if and only 
  \[
       \ord_x\omega\geq -1 \qquad\text{and} \qquad
        \Res_x\omega \in \FF_p,
  \]
  for all $x\in X$.
\end{lem}

\proof Let $x_1,\ldots,x_r$ be the set of poles of $\omega$ and set
$a_i:=\Res_{x_i}\omega$. After a change of coordinates we may assume
that $x_i\neq\infty$; then the point $x_i\in X$ is defined by $z=z_i$, for some
$z_i\in k$. 

Now suppose that $\omega$ has at most simple poles and that
$a_i\in\FF_p^\times$, for all $i$. Choose a lift $A_i\in\ZZ$ of $a_i$. Then
\[
     \omega = \sum_{i=1}^r \frac{a_i\,\diff z}{z-z_i} = \frac{\diff g}{g},
\]
with
\[
        g:= \prod_{i=1}^r \,(z-z_i)^{A_i}.
\]
This shows one direction of the claimed equivalence. The other direction is
obvious. 
\Endproof 

{\bf Proof of Proposition \ref{prop1}:}
We fix integers $m$ and $h$, with $m\geq 1$, $m|(p-1)$, $0<h<p$ and
$h\equiv -1\pmod{m}$. Write $h=mr-1$. The condition $h<p$ ensures that there
exist elements $z_1,\ldots,z_r\in\FF_p^\times$ such that 
\[
      z_{i,j}:=\zeta_m^jz_i\in\FF_p^\times, \quad i=1,\ldots,r,\;j=0,\ldots,m-1,
\]
are pairwise distinct. Define
\[
     \omega:= \frac{\diff z}{\prod_i (z^m-z_i^m)}.
\]
We claim that $\omega$ is a good deformation datum of conductor $h$. 

By construction, we have $\sigma^*\omega=\zeta_m\omega$, where
$\sigma$ is the automorphism of $X$ with $\sigma^*z=\zeta_mz$.
Furthermore, $\omega$ has exactly $mr$ simple poles and no zeroes on
$\AA_k^1\subset X$. It follows that $\omega$ has a zero of order
$h=mr-1$ at $\infty$. Finally, the residues of $\omega$ all lie in
$\FF_p$. Therefore, $\omega$ is logarithmic by Lemma \ref{lem1}. This
proves the claim and finishes the proof of Proposition \ref{prop1}.
\Endproof

\section{Applications to the lifting problem} \label{sec2}

Despite its simple proof, Proposition \ref{prop1} has a nontrivial
application to the local lifting problem, as already mentioned in the
introduction. 

First some notation. Let $P=\gen{\tau}\cong\ZZ/p\ZZ$ be cyclic of
order $p$, with generator $\tau$, and $C=\gen{\sigma}\cong\ZZ/m\ZZ$
cyclic of order $m$, with generator $\sigma$. Let
$\chi:C\to(\ZZ/p\ZZ)^\times$ be a character and $G:=P\rtimes_\chi C$
the corresponding semi-direct product (such that
$\sigma\tau\sigma^{-1}=\tau^{\chi(\sigma)}$, for $\sigma\in C$).  Let
$\phi:G\inj\Aut_k(k[[x]])$ be a local $G$-action. The {\em conductor}
of $\phi$ is defined as the positive integer
\[
        h := \ord_x\left(\frac{\tau(x)}{x}-1\right)
\]
(compare with \eqref{conductoreq}). Standard arguments from the theory of
local fields give the following restrictions on the conductor $h$ (see
\cite{Dp}).  Firstly, $h$ is prime to $p$.  Secondly, the order of $h$ in
$\ZZ/m\ZZ$ is equal to the order of the character $\chi$ (equivalently,
$\chi=\chi_0^h$, where $\chi_0:C\to K^\times$ is a primitive character of
order $m$). Conversely, if $h$ is a positive integer satisfying these two
conditions then there exists a local $G$-action with conductor $h$.

However, not every local $G$-action lifts to characteristic zero. The
following theorem gives a necessary and sufficient condition.

\begin{thm} \label{LiftThm}
  Let $\phi:G\inj\Aut_k(k[[x]])$ be a local $G$-action, where
  $G=P\rtimes_\chi C$ is the semi-direct product of a cyclic group $P$
  of order $p$ and a cyclic group $C$ of order $m$, with $(m,p)=1$.
  Let $h$ be the conductor of $\phi$. Then $\phi$ lifts to
  characteristic zero if and only if the following two conditions
  hold.
  \begin{enumerate}
  \item 
    The character $\chi:C\to(\ZZ/p\ZZ)^\times$ is either trivial
    or injective.  (Equivalently, the conductor $h$ is either
    divisible by $m$ or prime to $m$.)
  \item
    If $\chi$ is injective, then the congruence $h\equiv -1 \pmod{m}$ holds.
  \end{enumerate}
\end{thm}

There are two cases in which Conditions (i) and (ii) are automatically
true. The first case is when the character $\chi$ is trivial, i.e.
when $m|h$. Then $G$ is a cyclic group of order $pm$, and Theorem
\ref{LiftThm} says that every local $G$-action lifts to characteristic
zero. We recover a result from \cite{GM2}.

The second case is when $m=2$ and $\chi$ is not trivial (or
equivalently: $h$ is odd). Then $G$ is the dihedral group of order
$2p$. Again, Condition (i) and (ii) are automatically verified, and
Theorem 1 says that every local $G$-action lifts to characteristic
zero. This is the main result of \cite{Dp}.

In all other cases (i.e.\ for $m>2$ and $m\nmid h$) the congruence
$h\equiv -1\pmod{m}$ is strictly stronger than the condition that $h$
is prime to $m$ (i.e.\ that $\chi$ is injective). This means that not
every local $G$-action lifts to characteristic zero.

\vspace{1ex}
{\bf Proof of Theorem \ref{LiftThm}:} Having available Proposition
\ref{prop1}, the theorem is a straightforward consequence of the
methods of \cite{Dp}. We only sketch the argument. 

By \cite{Dp}, Corollary 3.7, liftability of $\phi$ is equivalent to
the existence of a {\em Hurwitz tree of type $(C,\chi)$, conductor $h$
  and discriminant $\delta=0$}. Hence the necessity of Conditions
(i) and (ii) of Theorem \ref{LiftThm} follows from \cite{Dp}, Lemma
3.3.(v).

To prove sufficiency of Condition (i) and (ii), it is natural to
distinguish the case of small conductor ($h<p$) from the case of large
conductor ($h>p$).

Suppose first that $h<p$. By \cite{GM}, Theorem III.3.1, a Hurwitz
tree of conductor $h<p$ is {\em irreducible}, i.e.\ it is the same
thing as a good deformation datum of conductor $h$. Therefore, the
sufficiency of Conditions (i) and (ii) follows from Proposition
\ref{prop1} (and necessity follows directly from Proposition
\ref{prop0}).

If $h>p$, there may not exist a good deformation datum of conductor
$h$. However, one can easily construct a (not necessarily irreducible)
Hurwitz tree of type $(C,\chi)$, conductor $h$ and discriminant
$\delta=0$, by adapting the proof of \cite{Dp}, Theorem 4.1 (which
deals with the case $m=2$ and $h>p$ odd) to the more general situation
considered here. This completes the proof of Theorem \ref{LiftThm}.
\Endproof

\section{The type} \label{type}

In the rest of this paper, we only consider primitive deformation data. 

We fix $m\geq 1$ as in Section \ref{sec1}. Let $\omega=f(z)\,\diff z$
be a good and primitive deformation datum. After a change of
coordinate, we may assume that the unique zero $x$ of $\omega$ is
located at the point $z=\infty$ (see the proof of Proposition
\ref{prop0}). This assumption allows us to identify the points different from
$x$ with elements of the field $k$.

Choose a system of representatives $z_1,\ldots, z_r$ of the
$\sigma$-orbits of the set of poles of $\omega$. By Lemma \ref{lem1}
we have
\[
      a_i :=\Res_{z_i}\omega \in\FF_p^\times.
\]
The tuple $\a:=(a_1,\ldots,a_r)$ is called the {\em type} of $\omega$.

Condition (a) of Definition \ref{def1} implies that $\omega$ has poles in
$z_{i,j}:=\zeta_m^jz_i$ and that
\[
     a_{i,j}:=\Res_{z_{i,j}}\omega = \zeta_m^{-jc}a_i, \quad j=0,\ldots,m-1,
\]
for some $c\in\ZZ$, prime to $m$. It follows that the points $z_{i,j}$ are
pairwise distinct. In particular, if $m>1$ then $z_{i,j}\in k^\times$. 

By the residue theorem we have
\begin{equation} \label{eq2}
   \sum_{i,j} a_{i,j} = 
      \Big(\sum_{i=1}^r a_i\Big)\cdot\Big(\sum_{j=0}^{m-1}\zeta_m^j\Big)  = 0.
\end{equation}
If $m>1$ then the second sum is zero, and the condition \eqref{eq2} coming
from the residue theorem is empty. 

\begin{defn} \label{def2}
  Let $r\geq 2$ be an integer such that $h:=mr-1\not\equiv 0\pmod{p}$.  A {\em
    type} of length $r$ is an $r$-tuple $\a=(a_1,\ldots,a_r)$, with
  $a_i\in\FF_p^\times$. If $m=1$ then we demand that $\sum_ia_i=0$. Two types
  $\a=(a_i)$ and $\a'=(a_i')$ of length $r$ are equivalent if there exists
  $\pi\in S_r$ and $c_i\in(\ZZ/m\ZZ)^\times$ with
  $a_i'=\zeta_m^{c_i}a_{\pi(i)}$, for $i=1,\ldots,r$.
\end{defn}

So the type of a deformation datum $\omega$ is well defined, up to
equivalence. We can now formulate a refined version of Problem \ref{prob1}.

\begin{prob} \label{prob2} Let $\a=(a_1,\ldots,a_r)$ be a type. Does there
  exist a good deformation datum $\omega$ of type $\a$?
\end{prob}

Let $\a=(a_1,\ldots,a_r)$ be a type. Choose elements $z_1,\ldots,z_r\in
k^\times$ such that the elements $z_{i,j}:=\zeta_m^jz_i\in k^\times$ are
pairwise distinct, for $i=1,\ldots,r$ and $j=0,\ldots,m-1$. Set
\[
     \omega := \sum_{i=1}^r \;\frac{m\,a_i\,z_i^{m-1}\,\diff z}{z^m-z_i^m}.
\]
This is clearly a differential form with at most simple poles in the points
$z_{i,j}$ satisfying $\sigma^*\omega=\zeta_m\omega$ (Condition (a) of
Definition \ref{def1}). A short computation shows that
\[
     \Res_{z_{i,j}}\omega = \zeta_m^{-j}a_i\in\FF_p^\times.
\]
By Lemma \ref{lem1}, $\omega$ is logarithmic (Condition (b) of
Definition \ref{def1}). Therefore, $\omega$ is a deformation datum of type
$\a$. It is uniquely determined by the choice of the poles $z_i$. 

For which choice of $z_1,\ldots,z_r$ is $\omega$ a good deformation datum? One
can rewrite
\[
    \omega = \frac{q(z)\,\diff z}{\prod_i(z^m-z_i^m)},
\]
where $q(z)\in k[z]$ is a polynomial of degree at most $h=mr-1$ whose
coefficients are polynomials in the $z_i$. The deformation datum $\omega$ is
good if and only if $q(z)$ is a constant. Setting the nonconstant coefficients
of $q(z)$ to zero, we obtain a system of polynomial equations in the unknown
$z_i$ whose solutions correspond to good deformation data of type $\a$. So
far, most attempts to construct good deformation data were based on explicitly
solving these equations.

In the remaining sections we suggest another approach to Problem
\ref{prob2}.

\section{A necessary condition} \label{sec4}

In this section we assume $m=1$. We formulate a necessary condition on a given
type $\a$ for the existence of good deformation data of type $\a$ (Corollary
\ref{cor2}). We then state our main result (Theorem \ref{thm1}) which says
that this condition is sufficient for large prime $p$. The proof is given in
the following section.

Let $\omega$ be a good deformation datum of type
$\a=(a_1,\ldots,a_r)$. We may assume that the point $z=\infty$ is the
unique zero of $\omega$. Then
\[
  \omega = \sum_{i=1}^r \frac{a_i\,\diff z}{z-z_i},
\]
with $z_1,\ldots,z_r\in k$, see the proof of Lemma \ref{lem1}. Note
that $\sum_i a_i=0$ and that the conductor of $\omega$ is $h=r-1$, which is
prime to $p$. 

For $i=1,\ldots,r$, choose a lift $A_i\in \ZZ$ of $a_i$ such that
\[
      \sum_{i=1}^r A_i = 0.
\]
We call $\A=(A_1,\ldots,A_r)$ a {\em lift} of the type $\a$. Then the rational
function
\begin{equation} \label{eq4}
     g := \prod_{i=1}^r (z-z_i)^{A_i}
\end{equation}
has the property that $\omega=\diff g/g$. 

We consider $g$ as a morphism $g:\PP^1_k\to\PP^1_k$. Let
$I\subset\{1,\ldots,r\}$ be the subset of indices $i$ with $A_i>0$.
Note that the degree of $g$ is equal to
\begin{equation} \label{eq5}
      \sum_{i\in I} A_i = -\sum_{i\not\in I}A_i 
           = \sum_{i=1}^r\max(A_i,0).
\end{equation}

\begin{prop} \label{prop2}
  The morphism $g:\PP^1_k\to\PP^1_k$ is branched
  exactly above $0,1,\infty$. More precisely:
  \begin{enumerate}
  \item
    $g$ is ramified at the point $z_i$, with ramification index $|A_i|$. Note
    that $g(z_i)=0$ if $i\in I$ and $g(z_i)=\infty$ otherwise.
  \item $g$ is ramified at $\infty$, with ramification index
    $\geq\min(h,p)$. It is tamely ramified at $\infty$ if and only if the
    ramification index is equal to $h$. Note that $g(\infty)=1$.
  \item
    $g$ is unramified at all other points.
  \end{enumerate}
  In particular, if $h<p$ then $g$ is tamely ramified everywhere.
\end{prop}

\proof Claim (i) is clear from the definition of $g$. Since $\omega$ is good
and $g(\infty)=1$, the differential form $\diff g=g\omega$ has a zero at
$\infty$ of order $h-1=r-2$ and no zero or pole on
$\AA^1_k-\{z_1,\ldots,z_r\}$. Claim (iii) follows immediately; Claim (ii) as
well, but one has to use the fact that the conductor $h$ is prime to $p$.  
\Endproof

Given a lifted type $\A=(A_1,\ldots,A_r)$, we set
\begin{equation} \label{nAkAeq}
   n_\A:=\sum_{i=1}^r \max(A_i,0), \qquad k_\A := {\rm gcd}(A_1,\ldots,A_r).
\end{equation}

\begin{cor} \label{cor2}
  Let $\omega$ be a good deformation datum of type $\a=(a_1,\ldots,a_r)$ (with
  $m=1$). Then for every lift $\A=(A_1,\ldots,A_r)$ of $\a$ we have
  \[
       k_\A\cdot\min(r-1,p) \leq n_\A.
  \]
\end{cor}

\proof Let $\A=(A_1,\ldots,A_r)$ be a lift of $\a$ and $k:=k_\A$. Then
$\A':=(A_1/k,\ldots,A_r/k)$ is a lift of $\a':=(a_1/k,\ldots,a_r/k)$
with $n_{\A'}=n_\A/k$ and $k_{\A'}=1$. Moreover,
$\omega':=k^{-1}\omega$ is a good deformation datum of type $\a'$. We
may therefore assume that $k_\A=1$.

Let $g:\PP^1_k\to\PP^1_k$ be the rational function associated to
$\omega$ and $\A$ by \eqref{eq4}. The degree of $g$ is equal to $n_\A$
and, by Proposition \ref{prop2} (ii), the ramification index of $g$ at
$\infty$ is $\geq \min(r-1,p)$. The inequality $\min(r-1,p)\leq n_\A$
follows immediately.
\Endproof

\begin{exa} \label{exa1} 
  Let $r\geq 4$ be even and $p>r/2$ a prime. Let
  \[
        \a:=(\underbrace{1,\ldots,1}_{\text{$r/2$ times}},
             \underbrace{-1,\ldots,-1}_{\text{$r/2$ times}}).
  \]
  Then $k_\A=1$ and $n_\A=r/2<\min(p,r-1)$, for the obvious lift
  $\A$. So Corollary \ref{cor2} implies that
  there does not exist a good deformation datum of type $\a$.

  The special case $p=5$ and $r=4$ corresponds to \cite{Dp}, Lemma 4.2. This
  example had been suggested to us by D. Harbater in connection with the
  lifting problem.  

  For $p=3$ and $r=6$, there does exists a one-parameter family of
  good deformation data of type $\a$. This shows that the bound in
  Corollary \ref{cor2} is sharp, in some sense.
\end{exa}

\begin{exa} \label{exa2} Let $s\geq 3$ and $1\leq\alpha <s$ be integers and
  $p>s$ be a prime. Set $r:=s+2$ and
  \[
     \a :=(\alpha-s,\underbrace{1,\ldots,1}_\text{$s$ times},-\alpha).
  \]
  Then $k_\A=1$ and $n_\A=s<\min(r-1,p)$ for the obvious lift $\A$. Again,
  Corollary \ref{cor2} shows that there does not exist a good
  deformation datum of type $\a$.  This result is proved in \cite{GM},
  Example 4.4, by explicit computations.
\end{exa} 

Here is our main result. 

\begin{thm} \label{thm1} Let $\a=(a_1,\ldots,a_r)$ be a type (with
  $m=1$). Suppose that there exists a lift $\A$ of $\a$ such that 
  \[
       k_\A\cdot (r-1) \leq n_\A < k_\A\cdot p.
  \]
  Then there exists a good deformation datum of type $\a$.
\end{thm}

The proof of Theorem \ref{thm1} is given at the end of the following
section. It is based on the idea, explained in the introduction, to construct
good deformation data from the reduction modulo $p$ of suitable Belyi maps.

\section{Good deformation data from Belyi maps} \label{belyi}

A {\em Belyi morphism} is a finite morphism $f:Y\to\PP^1_{\QQb}$, where $Y$ is
a smooth projective curve defined over the field $\QQb$ of algebraic numbers,
such that $f$ is branched exactly over the three points $0$, $1$ and $\infty$.

Belyi morphisms can be described by some purely topological
data. For instance, one can describe them by drawing a {\em dessin
  d'enfants}. An equivalent description, more useful for computations, is the
following. 

\begin{defn}
  Fix an integer $n\geq 1$. A {\em generating system} of degree $n$ is a
  triple $\ssigma=(\sigma_1,\sigma_2,\sigma_3)$ of permutations $\sigma_i\in
  S_n$ with the following properties:
  \begin{itemize}
  \item
    $\sigma_1\sigma_2\sigma_3=1$, and
  \item
    $G:=\gen{\sigma_1,\sigma_2,\sigma_3}\subset S_n$ acts transitively on
    $\{1,\ldots,n\}$. 
  \end{itemize}
  Two generating systems $\ssigma=(\sigma_i)$ and $\ssigma'=(\sigma_i')$ are
  called equivalent if there is a permutation $\tau\in S_n$ such that
  $\sigma_i'=\tau\sigma_i\tau^{-1}$, for $i=1,2,3$.

  The {\em combinatorial type} of a generating system $\ssigma$ is the triple
  of conjugacy classes of the elements $\sigma_i$ in $S_n$:
  \[
      C(\ssigma)=(C(\sigma_1),C(\sigma_2),C(\sigma_3)).
  \]
\end{defn}

By a standard construction (see e.g.\ \cite{Schneps}), we can associate to a
generating system $\ssigma$ a Belyi morphism $f:Y\to\PP^1_{\QQb}$. This gives
a one-to-one correspondence between equivalence classes of generating systems
of degree $n$ and Belyi morphisms of degree $n$.

If $f$ is the Belyi morphism corresponding to the generating system $\ssigma$,
then the combinatorial type $C(\ssigma)$ can be read off from the ramification
structure of $f$, as follows. Write $(x_1,x_2,x_3):=(0,1,\infty)$. Let
$f^{-1}(x_i)=\{y_1,\ldots,y_{s_i}\}$ be the fiber above $x_i$, and let $e_j$
denote the ramification index of $f$ at the point $y_j$. Then $C(\sigma_i)$ is
the conjugacy class of permutations with cycle structure $(e_1,\ldots,e_{s_i})$.
Note that we can compute the genus of the curve $Y$ solely from the
combinatorial type, using the Riemann--Hurwitz formula.

We say that a triple $C=(C_1,C_2,C_3)$ of conjugacy classes
of $S_n$ is {\em realizable} if there exists a generating system $\ssigma$
with $C=C(\ssigma)$. This means that there exists a Belyi morphism of degree
$n$ with combinatorial type $C$.

\vspace{1ex}
Let $\A=(A_1,\ldots,A_r)$ be a lifted type. Let
$C(\A)=(C_1,C_2,C_3)$ denote the following triple of conjugacy classes
of the symmetric group $S_{n_\A}$:
\begin{equation} \label{eq6}
    C_1 :=(A_i \mid i\in I),\quad C_2 :=(r-1), \quad
    C_3:=(-A_i \mid i\not\in I)
\end{equation}
(a conjugacy class in $S_n$ is given as a list of natural numbers indicating
the length of the cycles of a permutation, and in which $1$s can be omitted).
Note that the genus of $C(\A)$ is zero.

\begin{lem} \label{belyilem} 
  The triple of conjugacy classes $C(\A)$ is
  realizable if and only if the inequality
  \begin{equation} \label{belyilemeq1}
      k_\A\cdot(r-1) \leq n_\A.
  \end{equation}
  holds. Here $n_\A$ and $k_\A$ are defined by \eqref{nAkAeq}. 
\end{lem}

\proof 
It is no restriction to assume that
\begin{equation} \label{belyilemeq2}
   A_1,\ldots,A_s>0, \qquad A_{s+1},\ldots,A_r<0,
\end{equation}
for some $1\leq s<r$. 

Suppose that $f$ is a Belyi map of combinatorial type $C(\A)$. We may assume
that $f(\infty)=1$ and that $\infty$ is the unique ramification point above
$1$. We consider the inverse image $f^{-1}([-\infty,0])\subset\CC$ as a {\em
  dessins d'enfants}, i.e.\ a connected bicolored plane graph, as follows:
$f^{-1}(0)\subset\CC$ is the set of black vertices, $f^{-1}(\infty)$ is the
set of white vertices and the connected components of $f^{-1}((-\infty,0))$
are the edges of the graph. The valency list is $(A_1,\ldots,A_s)$ for the
black and $(-A_{s+1},\ldots,-A_r)$ for the white vertices.  The connected
components of $\CC-f^{-1}([-\infty,0])$ are in bijection with the fiber
$f^{-1}(1)$. In particular, the (unique) unbounded component corresponds to
$\infty$.

The fact that $f$ is unramified above $1$ except in the point $\infty$ (where
it is ramified of order $r-1$) is equivalent to the following: ($*$) the
boundary of every bounded connected component of $\CC-f^{-1}([-\infty,0])$
consists exactly of one black, one white vertex and two edges.

This condition is equivalent to the condition ($*$'): if we remove, for all
pairs of vertices of different colors, all edges joining the two vertices
except for one, we obtain a plane tree. We can thus associate to the Belyi map
$f$ a {\em bicolored weighted plane tree} $T_f$. Here the weight associated to
a edge joining the pair of vertices $(v_1,v_2)$ is the number of edges joining
$(v_1,v_2)$ in the original graph $f^{-1}([-\infty,0])$.

\begin{figure} \label{Atree}
\begin{center}

\setlength{\unitlength}{1.2mm}

\begin{picture}(80,20)

\put(0,5){\circle*{1}}
\put(10,5){\circle{1}}
\put(20,5){\circle*{1}}
\put(20,15){\circle{1}}
\put(30,5){\circle{1}}

\qbezier(0,5)(5,10)(10,5)
\qbezier(0,5)(5,0)(10,5)
\put(10,5){\line(1,0){10}}
\put(20,5){\line(1,0){10}}
\qbezier(20,5)(15,10)(20,15)
\qbezier(20,5)(25,10)(20,15)

\put(50,5){\circle*{1}}
\put(60,5){\circle{1}}
\put(70,5){\circle*{1}}
\put(70,15){\circle{1}}
\put(80,5){\circle{1}}

\put(50,5){\line(1,0){10}}
\put(60,5){\line(1,0){10}}
\put(70,5){\line(1,0){10}}
\put(70,5){\line(0,1){10}}

\put(54.3,6){$2$}
\put(64.3,6){$1$}
\put(74.3,6){$1$}
\put(68,9.5){$2$}

\end{picture}

\end{center}
\caption{A dessins of type $(2,4,-3,-2,-1)$}
\end{figure}
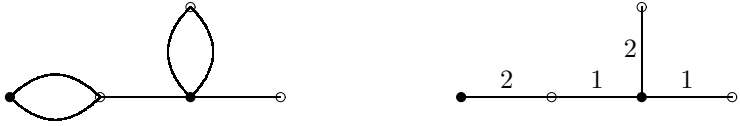

Conversely, given a bicolored weighted plane tree $T$ (where the weights of
the edges are positive integers), there exists a Belyi map $f$ of type $\A$,
unique up to isomorphism, such that $T_f$ is homeomorphic to $T$, see e.g.\
\cite{Schneps}. Here the tuple $\A=(A_1,\ldots,A_r)$ is the `valency list' of
$T$ (black vertices give positive numbers, white vertices give negative
numbers), and the valency of a vertex $v$ is the sum of the weights of all
edges adjacent to $v$. See Figure \ref{Atree} for a tree with valency list
$(2,4,-3,-2,-1)$ and the corresponding dessin. 

So far, we have translated the lemma into the following

\vspace{1ex} {\bf Claim:} There exists a bicolored weighted plane tree $T$
with valency list $\A$ if and only if the inequality \eqref{belyilemeq1}
holds.

\vspace{1ex} We will only prove the `if' part of the claim and leave the `only
if' part as an exercise. So we fix a lifted type $\A=(A_1,\ldots,A_r)$ such
that \eqref{belyilemeq1} and \eqref{belyilemeq2} hold and try to construct a
tree $T$ with valency list $\A$ (we abbreviate `bicolored weighted plane
tree' to `tree').

To begin, we observe that if $T$ is a tree with valency list $\A$, then the
gcd of the valencies, denoted $k_\A$ above, also divides the weight of every
edge of $T$ (proof by induction on the number of vertices). Multiplying
all the weights of a tree by a constant factor has the effect of multiplying
all valencies by the same factor. We may therefore assume that $k_\A=1$.

We now prove the claim by induction on the number $\min(s,r-s)$. Suppose first
that $\min(s,r-s)=1$, i.e.\ $s\in\{1,r-1\}$. Then one can draw a `star shaped'
tree $T$ with valency list $\A$, as follows. Say $s=1$. We choose one black
vertex $v_1\in\CC$ and $r-1$ white vertices $v_2,\ldots,v_r\in\CC$, positioned
in counterclockwise order on a circle around $v_1$. For $i=2,\ldots,r$ we then
join $v_1$ with $v_i$ by a straight line, which we regard as an edge of weight
$-A_i$.

So we may assume that $s,r-s\geq 2$. Then the $A_i$ cannot have all the same
absolute value. Therefore, after reordering the $A_i$ (and possibly changing
signs) we may also assume
\begin{equation} \label{belyilemeq3}
   A_1\leq A_2,\ldots,A_s,-A_{s+1},\ldots,-A_{r-1}, \quad A_1<-A_r.
\end{equation}
Set
\[
    m_i := {\rm gcd}(A_2,\ldots,A_i+A_1,\ldots,A_r),
    \quad\text{for $i=s+1,\ldots,r$.}
\]
The assumption $k_\A={\gcd}(A_1,\ldots,A_r)=1$ implies, by an easy argument,
that the numbers $m_{s+1},\ldots,m_r$ are pairwise relatively
prime. Moreover, if $A_i=-A_1$ for one index $i$, then $m_j=1$ for all
$j\neq i$. We conclude
\begin{equation} \label{belyilemeq4}
   \big(\prod_{j=s+1}^r m_j\big) \mid A_i, \quad
     \text{for $i=2,\ldots,s$}
\end{equation}
and 
\begin{equation} \label{belyilemeq5}
   \big(\prod_{j=s+1 \atop j\neq i}^r m_j\big) \mid A_i, \quad
      \text{for $i=s+1,\ldots,r$}.
\end{equation}
Furthermore, we can choose an index $i_0\in\{s+1,\ldots,r\}$ such that
\[
       m_{i_0}=m:=\min(m_{s+1},\ldots,m_r) \quad\text{and}\quad
        A_{i_0}+A_1<0.
\]

Set
\[
     \A':=(A_2,\ldots,A_{i_0}+A_1,\ldots,A_r)\in\ZZ^{r-1}.
\]
By construction, this is a lifted type with $r'=r-1<r$ entries and
\[
     n_{\A'} = n_{\A}-A_1, \qquad k_{\A'}=m
\]
(with the obvious notation). Then \eqref{belyilemeq4} implies that
\begin{equation} \label{belyilemeq6}
  n_{\A'} = A_2+\ldots+A_s \geq (s-1)\,m^{r-s},
\end{equation}
while \eqref{belyilemeq5} implies
\begin{equation} \label{belyilemeq7}
  n_{\A'} = -\big(\sum_{i=s+1 \atop i\neq i_0}^r A_i\big) -(A_{i_0}+A_1)
        \geq (r-s-1)\,m^{r-s-1}.
\end{equation}
We claim that
\begin{equation} \label{belyilemeq8}
   k_{\A'}\cdot(r'-1) \leq n_{\A'}.
\end{equation}
 
To prove the inequality \eqref{belyilemeq8} we suppose that it is false, i.e.\
that $m\,(r-2)>n_{\A'}$.  Then \eqref{belyilemeq6} implies
\begin{equation} \label{belyilemeq9}
  m^{r-s-1} < \frac{r-2}{s-1} = 1+\frac{r-s-1}{s-1}
\end{equation}
and \eqref{belyilemeq7} implies
\begin{equation} \label{belyilemeq10}
  m^{r-s-2} < \frac{r-2}{r-s-1} = 1+\frac{s-1}{r-s-1}
\end{equation}
Using the assumption $s,r-s\geq 2$, it is easy to deduce from
\eqref{belyilemeq9} and \eqref{belyilemeq10} that $m=1$.  But then the assumed
falsity of \eqref{belyilemeq8} means that $r-2>n_\A-A_1$, which implies
$k_{\A}\,(r-1)=r-1>n_\A$. This is a contradiction to \eqref{belyilemeq1} and
shows that \eqref{belyilemeq8} is true.\footnote{Although discovered
  independently, our argument to show \eqref{belyilemeq8} is almost identical
  to the argument given in \cite{Zannier}, proof of Theorem 1.}

It follows from \eqref{belyilemeq8} and the induction hypothesis that there
exists a tree $T'$ with valency list $\A'$. Let $v$ be a white vertex of
$T'$ with valency $-(A_{i_0}+A_1)$. Let $v'$ be any point in $\CC$ not lying
on $T'$. Let $T$ denote the tree obtained by adding to $T'$ the point $v'$ as a
black vertex and an edge with weight $A_1$ joining $v$ and $v'$. Then $T$
has valency list $\A$. This completes the proof of the lemma.
\Endproof

{\bf Proof of Theorem \ref{thm1}:} Let
$\a=(a_1,\ldots,a_r)\in(\FF_p^\times)^r$ be a type. Suppose that there
exists a lift $\A=(A_1,\ldots,A_r)$ such that
\[
       k_\A\cdot (r-1) \leq n_\A < k_\A\cdot p.
\]
We have to show that there exists a good deformation datum $\omega$ of
type $\a$. Replacing $\a$ by $\a':=\a/k_\A$, we may assume that
$k_\A=1$ (compare with the proof of Corollary \ref{cor2}). 

By Lemma \ref{belyilem}, there exists a Belyi map
$f:\PP^1_\QQb\to\PP^1_\QQb$ of combinatorial type $C(\A)$. After a
suitable choice of coordinates, we may assume that $f(\infty)=1$ and
that $\infty$ is the unique ramification point above $1$ (with
ramification index $r-1$). Now we have
\[
      f = \prod_{i=1}^r (z-z_i)^{A_i},
\]
for certain pairwise distinct algebraic numbers
$z_1,\ldots,z_r\in\QQb$. 

Let $v$ be a place of $\QQb$ above the prime $p$. Let $k$ denote the residue
field of $v$. By our hypothesis, the degree of $f$ is $n_\A<p$. Therefore, $f$
has {\em good reduction} at $v$, by Grothendieck's theory of the tame
fundamental group. In our concrete situation, this means that the algebraic
numbers $z_i$ are $v$-integers and that, moreover, the reduction
\[
      g = \prod_{i=1}^r (z-\bar{z}_i)^{A_i}\in k(z)
\]
of $f$ modulo $v$ is at most tamely ramified and branched exactly over $0$,
$1$, $\infty$, with ramification structure identical to that of $f$. It
follows that 
\[
     \omega:= \frac{\diff g}{g}
\]
is a good deformation datum of type $\a$. This proves the theorem.
\Endproof

\section{The case $m=2$}  \label{m=2}

Now suppose that $m=2$. Let $p$ be an odd prime and $\a=(a_1,\ldots,a_r)$ be a
type. A deformation datum of type $\a$ is of the form
\[
   \omega = \sum_{i=1}^r \Big(\frac{a_i\,\diff z}{z-z_i} 
                            - \frac{a_i\,\diff z}{z+z_i} \Big)
          = \sum_{i=1}^r \;\frac{2\,a_i\,z_i\,\diff z}{z^2-z_i^2},
\]
where $\pm z_i$ are pairwise distinct elements of $k^\times$. We suppose that
$\omega$ is good, i.e.\ that it has a single zero of order $h:=2r-1$ at
$\infty$. Moreover, we assume that $p>h$. 

Choose a lift $\A=(A_1,\ldots,A_r)$ of $\a_i$ with $A_i>0$, and set
\[
    g:=\prod_{i=1}^r \Big(\frac{z-z_i}{z+z_i}\Big)^{A_i}.
\]
Then $\omega=\diff g/g$. Furthermore, we have
\[
     g(-z) = \prod_{i=1}^r \Big(\frac{z+z_i}{z-z_i}\Big)^{A_i}
           = \frac{1}{g(z)}.
\]
Therefore, there exists a rational function $\tilde{g}\in k(x)$ such that
\[
    \Big(\frac{g(-z)-1}{g(-z)+1}\Big)^2 =
    \Big(\frac{g(z)-1}{g(z)+1}\Big)^2 = \tilde{g}(z^2).
\]

\begin{prop} \label{prop4} 
  The morphism $\tilde{g}:\PP^1_k\to\PP^1_k$ has degree 
  \[
       n=n_\A := \sum_{i=1}^r A_i.
  \]     
  It is at most tamely ramified and branched exactly above $0,1,\infty$. More
  precisely:
  \begin{enumerate}
  \item The fiber $\tilde{g}^{-1}(1)$ consists of the points $x=z_i^2$, for
    $i=1,\ldots,r$. The ramification index of $\tilde{g}$ at the point
    $x=z_i^2$ is equal to $A_i$. 
  \item 
    Suppose that $n$ is even. Then the fiber $\tilde{g}^{-1}(0)$ consists of
    the point $x=\infty$, with ramification index $h=2r-1$, of the point $x=0$,
    with ramification index $1$, and of $(n-2r)/2$ points with ramification
    index $2$. The fiber $\tilde{g}^{-1}(\infty)$ consists of $n/2$ points
    with ramification index $2$. 
  \item Suppose that $n$ is odd. Then the fiber $\tilde{g}^{-1}(0)$ consists
    of the point $x=\infty$, with ramification index $h=2r-1$ and of
    $(n-2r+1)/2$ points with ramification index $2$. The fiber
    $\tilde{g}^{-1}(\infty)$ consists of the point $x=0$, which is unramified,
    and of $(n-1)/2$ points with ramification index $2$. 
  \end{enumerate}
\end{prop}
  
\proof
We have a commutative diagram
\[\begin{CD}
   \PP^1_z @>{\psi}>> \PP^1_x  \\
   @V{g}VV        @VV{\tilde{g}}V \\
   \PP^1_v @>{\phi}>> \PP^1_u,\\
\end{CD}\]
where $\psi(z)=z^2=:x$ and $\phi(v)=((v-1)/(v+1))^2=:u$. In particular,
$\phi(0)=\phi(\infty)=1$, $\phi(1)=0$ and $\phi(-1)=\infty$. So 
\[
   \tilde{g}^{-1}(1) = \psi(g^{-1}(0)\cup g^{-1}(\infty))=\{z_1^2,\ldots,z_r^2\}.
\]
Moreover, since $g$ is ramified at $z=\pm z_i$ of order $A_i$ (Lemma \ref{lem1})
and $\psi$ is unramified at $z=\pm z_i$, we see that the ramification index of
$\tilde{g}^{-1}$ at $x=z_i^2$ is equal to $A_i$. This proves (i). 

One proves (ii) and (iii) in the same manner. The case distinction comes from
\[
      g(0) = (-1)^n \quad\Rightarrow\quad
        \tilde{g}(0) = \begin{cases} \;\;0, & \text{if $n$ is even,}\\
                                     \;\;\infty, & \text{if $n$ is odd.}
                       \end{cases}
\]
\Endproof

For a lift $\A=(A_1,\ldots,A_r)$ of $\a$ we set $C(\A)=(C_1,C_2,C_3)$, where
\[
  C_1 := (2r-1,\underbrace{2,\ldots,2}_{[(n-2r+1)/2]}), \quad 
  C_2 := (A_1,\ldots,A_r), \quad
  C_3 := (\underbrace{2,\ldots,2}_{[n/2]})
\]
are conjugacy classes in $S_{n_\A}$.  With this notation, we have obtain the
following version of Corollary \ref{cor2} and Theorem \ref{thm1} for $m=2$.

\begin{thm} \label{thm2}
    Let $\a$ be a type, with $m=2$ and $p>h=2r-1$. 
  \begin{enumerate}
  \item If there exists a good deformation datum $\omega$ of type $\a$, then
    the triple of conjugacy classes $C(\A)$ is
    realizable, for every lift $\A$ of $\a$.
  \item
    Suppose that there exists a lift $\A$ of $\a$ such that
    \begin{enumerate}
    \item
      $C(\A)$ is realizable, and
    \item
      $n_\A<p$.
    \end{enumerate}
    Then there exists a good deformation datum $\omega$ of type $\a$.
  \end{enumerate}
\end{thm} 

We have not tried to find an explicit criterion for realizability in the above
situation (such as Lemma \ref{belyilem}). Instead, we give an example:

\begin{exa} \label{exa4} For an arbitrary $r\geq 1$, set $m:=2$ and
  \[
       \A:=(\underbrace{1,\ldots,1}_{r-1},r). 
  \]
  Then $n=h=2r-1$ and the triple $C(\A)$ consists of the following conjugacy
  classes in $S_n$: $C_1$ is the class of all $n$-cycles, $C_2$ the class of
  $r$-cycles and $C_3$ the class of the product of $r-1$ disjoint $2$-cycles.
\end{exa}
  
\begin{lem} \label{lem4}
  Up to equivalence, there is a unique generating system
  $\ssigma$ of combinatorial type $C(\A)$.
\end{lem}

\proof
The existence of $\ssigma$ follows from the identity
\[
   (1\;2\,\ldots\, n) \cdot (n\;\,n-2\,\ldots\,3\;1) 
     = (2\;3)(4\;5)\cdots(n-1\;\,n).
\]
The proof of the uniqueness is an easy exercise. 
\Endproof

One can also convince oneself of the truth of Lemma \ref{lem4} by drawing the
corresponding dessin d'enfant. If $f:\PP^1\to\PP^1$ is a Belyi map of
combinatorial type $C(\A)$, with $\A$ as above, then the inverse image
$T:=f^{-1}([1,\infty])\subset\CC$ is a {\em pre-clean plane tree} with valency
list $(h,1,\ldots,1)$ (with $r$ repetitions of $1$), see e.g.\
\cite{Schneps}. It is clear that there is, up to isomorphism, a unique
realization of such a tree:
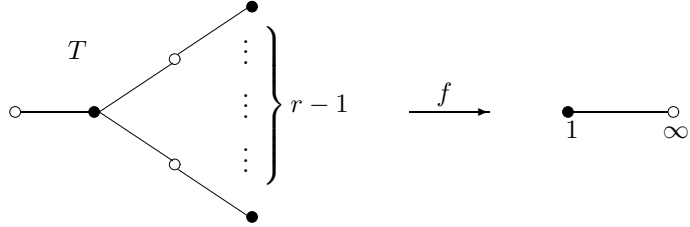
\begin{center}

\setlength{\unitlength}{0.7mm}

\begin{picture}(140,40)

\put(10,20){\circle{2}}
\put(11,20){\line(1,0){13}}
\put(25,20){\circle*{2}}
\put(26,20){\line(3,2){14}}
\put(40.3,30){\circle{2}}
\put(26,20){\line(3,-2){14}}
\put(40.3,10){\circle{2}}
\put(41,31){\line(3,2){13.5}}
\put(55,40){\circle*{2}}
\put(41,9){\line(3,-2){13.5}}
\put(55,0){\circle*{2}}

\put(20,30){$T$}

\put(50,20){$\left.\begin{array}{c} \vdots \\ \vdots \\ 
                                    \vdots\end{array}\right\} r-1$}

\put(85,20){\vector(1,0){15}}
\put(90,22){$f$}

\put(115,20){\circle*{2}}
\put(116,20){\line(1,0){18}}
\put(135,20){\circle{2}}
\put(114.5,15){$1$}
\put(133,15){$\infty$}

\end{picture}

\end{center}
       
\begin{cor}
  Let $r\geq 1$ and $p>n:=2r-1$ a prime number. Set $m:=2$. Then there exists
  a unique good deformation datum of type $\a:=(1,\ldots,1,r)$, defined over
  $\FF_p$. 
\end{cor}

This is essentially what is proved in \cite{Dp}, Section 5.3, by showing that
a certain system of equations has a unique solution. The above proof is
clearly shorter and more elegant.

\vspace{1.5cm}\noindent\small
irene.bouw@uni-ulm.de\\
wewers@math.uni-hannover.de\\
zapponi@math.jussieu.fr

\end{document}